\documentclass[12pt,reqno]{article}

\usepackage[english]{babel}
\usepackage[latin1]{inputenc}
\usepackage[dvips,final]{graphics}
\usepackage{amsmath,amsfonts,amssymb,amsthm,amscd,array,
stmaryrd,mathrsfs,mathdots,epigraph}
\usepackage[makeroom]{cancel}
\usepackage[all]{xy}
\usepackage{url}
\usepackage{pythonhighlight}
\usepackage{pgfplots}
\pgfplotsset{compat=1.15}
\usepackage{mathrsfs}
\usetikzlibrary{arrows}
\definecolor{qqqqtt}{rgb}{0,0,0.2}
\definecolor{ffqqqq}{rgb}{1,0,0}
\usepackage{blindtext}
\usepackage{textcomp}
\usepackage[final]{epsfig}
\usepackage{color}
\usepackage{mathabx}

\def\a{\alpha}

\def\e{\varepsilon}
\def\g{\gamma}

\def\s{\sigma}

\theoremstyle{plain}

\newtheorem{thm}{Theorem}
\newtheorem{conj}{Conjecture}
\newtheorem{lem}[thm]{Lemma}

\newcommand{\R}{\mathbb{R}}
\newcommand{\Z}{\mathbb{Z}}

\newcommand{\N}{\mathbb{N}}

\newcommand{\Q}{\mathbb{Q}}

\newcommand{\QP}{{\mathbb{QP}}}

\usepackage{amssymb}
\usepackage{amsmath}
\usepackage{amsthm}
\usepackage{amsfonts}
\usepackage{amscd}
\usepackage{graphicx}

\usepackage[colorlinks=true,
linkcolor=webgreen,
filecolor=webbrown,
citecolor=webgreen]{hyperref}

\definecolor{webgreen}{rgb}{0,.5,0}
\definecolor{webbrown}{rgb}{.6,0,0}

\usepackage{color}
\usepackage{fullpage}
\usepackage{float}

\usepackage{graphics}
\usepackage{latexsym}
\usepackage{epsf}
\usepackage{breakurl}

\newcommand{\seqnum}[1]{\href{https://oeis.org/#1}{\rm \underline{#1}}}

\begin{document}

\theoremstyle{plain}
\newtheorem{theorem}{Theorem}
\newtheorem{corollary}[theorem]{Corollary}
\newtheorem{lemma}[theorem]{Lemma}
\newtheorem{proposition}[theorem]{Proposition}

\theoremstyle{definition}
\newtheorem{definition}[theorem]{Definition}
\newtheorem{example}[theorem]{Example}
\newtheorem{conjecture}[theorem]{Conjecture}

\theoremstyle{remark}
\newtheorem{remark}[theorem]{Remark}

\begin{center}
\vskip 1cm{\LARGE\bf Conditions of Positivity on a Shadow Markoff Tree
}
\vskip 1cm
\large
Nathan Bonin\\
Laboratoire de Math\'ematiques de Reims\\
U.F.R. Sciences Exactes et Naturelles\\
Universit\'e de Reims Champagne-Ardenne\\
Moulin de la Housse - BP 1039,
51687 Reims cedex 2\\
France\\
\href{mailto:nathan.bonin1@etudiant.univ-reims.fr}{\tt nathan.bonin1@etudiant.univ-reims.fr} \\
\end{center}

\vskip .2 in

\begin{abstract}
I study an analog of the Markoff equation introduced in a previous paper with Valentin Ovsienko and
formulate a conjecture about the necessary and sufficient conditions for the positivity of solutions to this equation.
This conjecture is based on computer experiments and a theorem that gives a partial answer.
\end{abstract}

\thispagestyle{empty}

\section{Introduction}

The classical Markoff equation~\cite{Mar} is the Diophantine equation
\begin{equation}
\label{MarEq}
a^2+b^2+c^2=3abc.
\end{equation}
{\it Markoff numbers} are elements of the triple $(a,b,c)$ of positive integer solutions to~\eqref{MarEq}. 
They were introduced by Andrey Markoff in 1879 and raised much interest in many branches of mathematics, such as 
number theory, topology, combinatorics, and mathematical physics; see, e.g.,~\cite{Bom,Rud,Pro,Zag}.

It is an easy part of the Markoff theorem, that 
all positive integer solutions to~\eqref{MarEq} can be obtained from the ``initial'' triple $(1,1,1)$ 
via a process called mutation, combined with permutations of $a,b,c$ (see, e.g.,~\cite{Aig}).
Given a triple $(a,b,c)$; the {\it mutation} with respect to $a$ is another triple of solutions $(a',b,c)$ with 
\begin{equation}
\label{muta}
a'=\dfrac{b^2+c^2}{a}.
\end{equation}
An equivalent expression is $a'=3bc-a$.
The positive triple of solutions to~\eqref{MarEq} are organized in the form of a tree called the {\it Markoff tree}.
Note that the term ``mutation'' is due to the relation to cluster algebra; see~\cite{FZ1,Pro}.
Note also that mutations are involutions, i.e., double mutation at $a$ is identity: $a''=a$.

The following analog of the Markoff equation was introduced by Bonin and Ovsienko~\cite{OB}:
\begin{equation}
\label{SMarEq}
A^2+B^2+C^2=\left(3-\sigma\e\right)ABC
\end{equation}
where $A,B$ and $C$ are called {\it dual numbers}, i.e.,
\begin{equation}
\label{ArbTrip}
A=a+\a\e,
\qquad
B=b+\beta\e,
\qquad
C=c+\g\e.
\end{equation}
Dual numbers are elements of a commutative algebra of the form $A=a+\a\e $
where $a,\a\in\R$ and $\e$ is a formal parameter such that $\e^2=0$.
Following the terminology of~\cite{Ovs}, Eqn.~\eqref{SMarEq} was called the Shadow Markoff equation by Bonin and Ovsienko~\cite{OB}.

The notion of ``shadow'' sequences of integers appeared in~\cite{Ovs} (also see~\cite{CO, Hon,Ves}), and was tested on
the sequence of Markoff numbers.
Every Markoff number $a$ is accompanied by another integer, $\a$, 
called the shadow of $a$. This process implies the choice of initial conditions. For a definition, see~\cite{Ovs}.

It was proved in~\cite{OB} that~\eqref{SMarEq} is the unique 
$\e$-deformed Markov equation of the form
$$
A^2+B^2+C^2=3ABC+P(A,B,C)\e,
$$
where $P$ is an arbitrary polynomial in $A,B,C$, which is stable under the {\it mutations} 
of dual numbers given by the same formula
$\mu_A: (A,B,C)\mapsto(A',B,C)$, where
\begin{equation}
\label{MutA}
A'=\dfrac{B^2+C^2}{A}.
\end{equation}
More explicitly, the mutation reads
\begin{equation}
\label{SMarmutBis}
a'=\frac{b^2+c^2}{a},
\qquad\qquad
\a'=\frac{-a'\a+2b\,\beta+2c\,\g}{a}.
\end{equation}
If $(A,B,C)$ is a solution to~\eqref{SMarEq}, then $(A',B,C)$ is also a solution.

Note that the integrality of the solutions after mutations is guaranteed by the Laurent phenomenon of~\cite{FZ} (also see~\cite{OZ}).
More general mutation rules in the situation with nilpotent parameters can be found in~\cite{OZ}.

It follows from the Markoff theorem (for details; see~
\cite[\cite{OB}, Proposition 2(i)]{Alces}
\cite{OB}, Proposition 2(i)) that
every triple~\eqref{ArbTrip} of integer solutions to~\eqref{SMarEq}, i.e.,
such that $a,b,c,\a,\beta,\g\in\Z$, can be obtained by a sequence of mutations~\eqref{MutA}
and permutations from the ``initial triple'' of the form
\begin{equation}
\label{InitTrip}
(A,B,C)=
(1+\alpha\e,\; 1+\beta\e, \; 1+\gamma\e),
\end{equation}
where $(\a,\beta,\g)$ are some integers.
Note that the constant $\s\in\R$ in~\eqref{SMarEq} is equal to the sum of the nilpotent parts of the initial values of $A,B,C$:
$$
\sigma=\alpha+\beta+\gamma.
$$
This follows directly from~\eqref{SMarEq} and~\eqref{InitTrip}; also see Bonin and Ovsienko ~\cite[p.\ 1487]{OB}

Every initial triple~\eqref{InitTrip} corresponds to a point $(\a:\beta:\g)$ in 
the rational projective plane~$\QP^2$.
This correspondence is natural since the mutations~\eqref{SMarmutBis} are linear on the nilpotent part of solutions.
For more details; see Section~\ref{ProofSec}.

Along the Markoff tree, $a,b,c$ remain positive while $\a,\beta,\g$ may become negative.
The problem of characterization of {\it positive integer solutions} to~\eqref{SMarEq} was formulated in~\cite{OB}.
The problem is to describe all the solutions to~\eqref{SMarEq} with
\begin{equation}
\label{PosCon}
a,b,c,\a,\beta,\g\in\N.
\end{equation}
It can be reformulated in terms of the
 initial triple~\eqref{InitTrip}
 such that under every series of mutations along the Markoff tree
the  triple $(\a,\beta,\g)$ remain positive.
In this paper, I study this problem and formulate a conjecture that gives a complete answer to it.

\begin{conj}
\label{MainConj}
Every triple $(A,B,C)$ of positive integer solutions to~\eqref{SMarEq} can be obtained
from the initial triple~\eqref{InitTrip}
where $(\a,\beta,\g)$ are positive integers such that $(\a:\beta:\g)\in\QP^2$ are the rational points that belong to the quadrilateral 
with vertices 
$$
\textstyle
\left\{(0:0:1); \quad (\dfrac{1}{2}:0:1)\equiv(1:0:2); \quad (1:1:1); \quad (0:2:1)\right\}
$$
depicted in the following figure
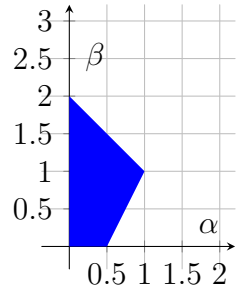
\begin{center}
\definecolor{zzttqq}{rgb}{0,0,1}
\begin{tikzpicture}[line cap=round,line join=round,>=triangle 45,x=1cm,y=1cm]
\begin{axis}[
x=1cm,y=1cm,
axis lines=middle,
ymajorgrids=true,
xmajorgrids=true,
ymin=-0.29868247168517575,
ymax=3.2156316474709685,
xmin=-0.36602351946208256,
xmax=2.204001897585711,
ytick={0,0.5,...,3},
xtick={0,0.5,...,2},]
\clip(-0.36602351946208256,-0.29868247168517575) rectangle (4.204001897585711,3.2156316474709685);
\fill[line width=2pt,color=zzttqq,fill=zzttqq,fill opacity=0.50000000149011612] (0,0) -- (0.5,0) -- (1,1) -- (0,2) -- cycle;
\draw (0.06610799040603345,2.8605427849359513) node[anchor=north west] {$\beta$};
\draw (1.5738817601573287,0.50538827311252061) node[anchor=north west] {$\alpha$};
\end{axis}
\end{tikzpicture}
\end{center}
via mutations \eqref{MutA} and permutations along the Markoff tree.
\end{conj}

In Section~\ref{ProofSec}, I will prove the following.

\begin{thm}
\label{TheThm}
Every triple leading to a positive integer tree of solutions to~\eqref{SMarEq} can be obtained from some initial triple~\eqref{InitTrip}
with the point $(\a:\beta:\g)\in\Q$P$^2$ that belongs to some convex polygon
by a series of mutations and permutations along the Markoff tree.
\end{thm}
The main ingredient of the proof of this statement is the fact that $\alpha'$ 
in~\eqref{SMarmutBis} depends linearly on $(\alpha,\beta,\gamma)$.

I will also provide ``numeric computer-assisted computations'' that give evidence that the above conjecture is indeed true.
That is, the convex polygon is precisely the quadrilateral from Conjecture~\ref{MainConj}.

\section{The classic Markoff tree}

The solutions of the classical Markoff equation can be visualized in the form of a tree. 
Every Markoff triple $(a,b,c)$ labels three regions in the plane: 
$$
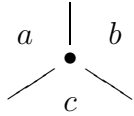
\xymatrix @!0 @R=0.3cm @C=0.3cm
{
&&&\ar@{-}[ddd]\\
\\
&a&&&&b\\
&&&\bullet\ar@{-}[llldd]\ar@{-}[rrrdd]\\
\\
&&&c&&&
}
$$
and the mutation of a triple $(a,b,c)$ corresponds to the following branchings:
$$
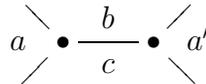
\xymatrix @!0 @R=0.3cm @C=0.3cm
{
\ar@{-}[rrdd]&&&&&&&&\\
&&&&b\\
a&&\bullet\ar@{-}[lldd]\ar@{-}[rrrr]&&&&
\bullet\ar@{-}[rruu]\ar@{-}[rrdd]&&a'\\
&&&&c\\
&&&&&&&&
}
$$

The classical Markoff tree:
$$
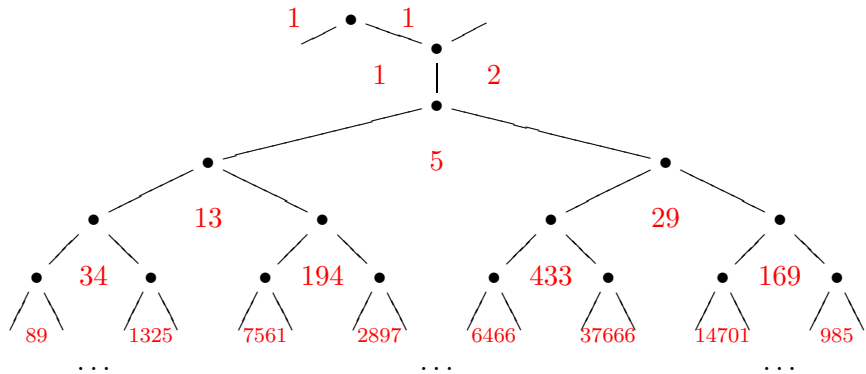
\begin{small}
\xymatrix @!0 @R=0.38cm @C=0.38cm
{&&&&&&&&&&&&\\
&&&&&&&&&&&{\textcolor{red}{1}}
&&\bullet\ar@{-}[lld]\ar@{-}[lu]&&{\textcolor{red}{1}}&&&\\
&&&&&&&&&&&&&&&&\bullet\ar@{-}[dd]\ar@{-}[rru]\ar@{-}[lllu]\\
&&&&&&&&&&&&&&{\textcolor{red}{1}}&&&&{\textcolor{red}{2}}\\
&&&&&&&&&&&&&&&&\bullet\ar@{-}[lllllllldd]\ar@{-}[rrrrrrrrdd]&&&&&&&&\\
&&&&&&&&&&&&&&&&\\
&&&&&&&&\bullet\ar@{-}[lllldd]\ar@{-}[rrrrdd]
&&&&&&&&{\textcolor{red}{5}}&
&&&&&&&\bullet\ar@{-}[lllldd]\ar@{-}[rrrrdd]&&&\\
&&&&&&&&
&&&&&&&&&&&&&&&&\\
&&&&\bullet\ar@{-}[lldd]\ar@{-}[rrdd]
&&&&{\textcolor{red}{13}}&&&&\bullet\ar@{-}[lldd]\ar@{-}[rrdd]
&&&&&&&&\bullet\ar@{-}[lldd]\ar@{-}[rrdd]
&&&&{\textcolor{red}{29}}&&&&\bullet\ar@{-}[lldd]\ar@{-}[rrdd]\\
&&&&&&&&&&&&&&&&&&&&&&&&&&&&\\
&&\bullet\ar@{-}[ldd]\ar@{-}[rdd]
&&{\textcolor{red}{34}}&&\bullet\ar@{-}[ldd]\ar@{-}[rdd]
&&&&\bullet\ar@{-}[ldd]\ar@{-}[rdd]
&&{\textcolor{red}{194}}&&\bullet\ar@{-}[ldd]\ar@{-}[rdd]
&&&&\bullet\ar@{-}[ldd]\ar@{-}[rdd]
&&{\textcolor{red}{433}}&&\bullet\ar@{-}[ldd]\ar@{-}[rdd]
&&&&\bullet\ar@{-}[ldd]\ar@{-}[rdd]
&&{\textcolor{red}{169}}&&\bullet\ar@{-}[ldd]\ar@{-}[rdd]\\
&&&&&&&&&&&&&&&&&&&&&&
&&&&&&&&&&&\\
&&{\textcolor{red}{\scriptstyle89}}
&&&&{\textcolor{red}{\scriptstyle1325}}
&&&&{\textcolor{red}{\scriptstyle7561}}
&&&&{\textcolor{red}{\scriptstyle2897}}
&&&&{\textcolor{red}{\scriptstyle6466}}
&&&&{\textcolor{red}{\scriptstyle37666}}
&&&&{\textcolor{red}{\scriptstyle14701}}
&&&&{\textcolor{red}{\scriptstyle985}}&&
\\
&&&&\ldots&&&&&&&&&&&&\ldots&&&&&&&&&&&&\ldots
}
\end{small}$$
is the standard infinite binary tree cutting the plane into regions labeled by all Markoff numbers.

The initial triple $(1,1,1)$ is the root of the tree.
It is followed by $(1,1,2)$ and $(1,2,5)$.
After that, any mutation is allowed.
The left branch consists of the Fibonacci numbers with odd indices  $F_{2k+1}$,
while the right branch is that of the odd Pell numbers  $P_{2k+1}$.

\section{Proof of the main theorem}\label{ProofSec}

In this section, I prove the main result, Theorem~\ref{TheThm}. It can be reformulated 
as follows.

\begin{thm}
\label{TheThmBis}
There exists a convex polygon $P$ in $\QP^2$ such that every solution to 
the Shadow Markoff equation~\eqref{SMarEq}
obtained from an initial triple~\eqref{InitTrip} is positive if and only if the point
$(\a:\beta:\g)$ belongs to $P$.
\end{thm}

\begin{proof}
Consider two initial triple, $T_1^0$ and~$T_2^0$ as in~\eqref{InitTrip}, 
and assume that both $T_1^0$ and~$T_2^0$ produce positive integer solutions to~\eqref{SMarEq}.
Take $(\lambda, \mu)=(\frac{p}{q},\frac{r}{s})$ positive rational numbers, such that  $\lambda+\mu=1$.
Consider the linear combination
$$
T^0=\lambda T_1^0 + \mu T_2^0.
$$
It is of the form
$$
T^0=
\Big(
1+\big(\dfrac{p}{q}\alpha_1+\dfrac{r}{s}\alpha_2\big)\e,\; 
1+\big(\dfrac{p}{q}\beta_1+\dfrac{r}{s}\beta_2\big)\e,\; 
1+\big(\dfrac{p}{q}\gamma_1+\dfrac{r}{s}\gamma_2\big)\e
\Big).
$$
The nilpotent part of $T^0$ corresponds to the following point in $\QP^2$:
$$
(ps\alpha+rq\alpha_2\,:\,ps\beta_1+rq\beta_2\,:\,ps\gamma_1+rq\gamma_2).
$$
To prove Theorem~\ref{TheThm}, one needs to show that the initial triple
$$
\tilde T^0:=
\big(
1+\left(ps\alpha_1+rq\alpha\right)\epsilon,\;
1+\left(ps\beta_1+rq\beta_2\right)\epsilon,\;
1+\left(ps\gamma_1+rq\gamma_2\right)\epsilon
\big),
$$
that can be viewed as the barycenter of $T_1^0$ and~$T_2^0$,
corresponds to positive integer solutions of~\eqref{SMarEq} after mutations along the Markoff tree.

The statement then follows from the linearity of the nilpotent part of the mutation~\eqref{SMarmutBis}.
Indeed, if $(a+\a\epsilon,b+\beta\epsilon,c+\gamma\epsilon)$ and 
$(a+\alpha'\epsilon,b+\beta'\e,c+\gamma'\epsilon)$ stay at the same place in 
the Markoff trees of $T_1^0$ and~$T_2^0$, respectively, then
at the same place in the Markoff tree of $\tilde T^0$ one gets the triple
$$
\big(a+\left(ps\a+rq\a'\right)\e,\;
b+\left(ps\beta+rq\beta'\right)\e,\;
c+\left(ps\g+rq\g'\right)\e\big).
$$
Hence the result.
\end{proof}

Theorem~\ref{TheThm} is proved.

\section{Four shadow Markoff trees}

Replacing the initial Markoff triple $(1,1,1)$ by an initial
triple~\eqref{InitTrip}, one obtains a node
with the root labeled as follows.
$$
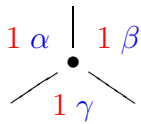
\xymatrix @!0 @R=0.3cm @C=0.3cm
{
&&&\ar@{-}[ddd]\\
\\
&\textcolor{red}{1}\;\textcolor{blue}{\alpha}&&&&\textcolor{red}{1}\;\textcolor{blue}{\beta}\\
&&&\bullet\ar@{-}[llldd]\ar@{-}[rrrdd]\\
\\
&&&\textcolor{red}{1}\;\textcolor{blue}{\gamma}&&&
}
$$
Then, following the mutations along the Markoff tree,
one obtains a tree of solutions to~\eqref{SMarEq}
where $\sigma=\alpha+\beta+\gamma$.
In this section, I present the trees corresponding to four vertices of the quadrilateral from Conjecture~\ref{MainConj}.

\subsection{The tree of the vertex $(0:0:1)$}
~\\
Taking the initial triple
$(1,\; 1, \; 1+\e)$, one has the following tree.
			$$
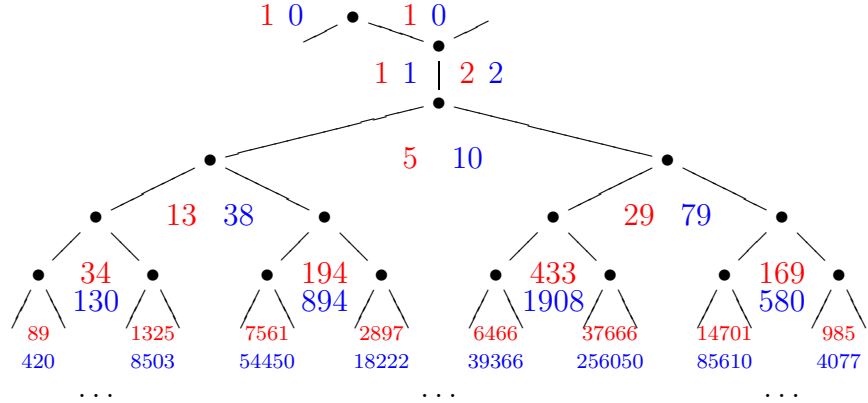
\xymatrix @!0 @R=0.38cm @C=0.38cm
{&&&&&&&&&&&&\\
&&&&&&&&&&{\textcolor{red}{1}}&{\textcolor{blue}{0}}
&&\bullet\ar@{-}[lld]\ar@{-}[lu]&&{\textcolor{red}{1}}&{\textcolor{blue}{0}}&&\\
&&&&&&&&&&&&&&&&\bullet\ar@{-}[dd]\ar@{-}[rru]\ar@{-}[lllu]\\
&&&&&&&&&&&&&&{\textcolor{red}{1}}&{\textcolor{blue}{1}}&&{\textcolor{red}{2}}&{\textcolor{blue}{2}}\\
&&&&&&&&&&&&&&&&\bullet\ar@{-}[lllllllldd]\ar@{-}[rrrrrrrrdd]&&&&&&&&\\
&&&&&&&&&&&&&&&&\\
&&&&&&&&\bullet\ar@{-}[lllldd]\ar@{-}[rrrrdd]
&&&&&&&{\textcolor{red}{5}}&&{\textcolor{blue}{10}}
&&&&&&&\bullet\ar@{-}[lllldd]\ar@{-}[rrrrdd]&&&\\
&&&&&&&&
&&&&&&&&&&&&&&&&\\
&&&&\bullet\ar@{-}[lldd]\ar@{-}[rrdd]
&&&{\textcolor{red}{13}}&&{\textcolor{blue}{38}}&&&\bullet\ar@{-}[lldd]\ar@{-}[rrdd]
&&&&&&&&\bullet\ar@{-}[lldd]\ar@{-}[rrdd]
&&&{\textcolor{red}{29}}&&{\textcolor{blue}{79}}&&&\bullet\ar@{-}[lldd]\ar@{-}[rrdd]\\
&&&&&&&&&&&&&&&&&&&&&&&&&&&&\\
&&\bullet\ar@{-}[ldd]\ar@{-}[rdd]
&&{\textcolor{red}{34}}&&\bullet\ar@{-}[ldd]\ar@{-}[rdd]
&&&&\bullet\ar@{-}[ldd]\ar@{-}[rdd]
&&{\textcolor{red}{194}}&&\bullet\ar@{-}[ldd]\ar@{-}[rdd]
&&&&\bullet\ar@{-}[ldd]\ar@{-}[rdd]
&&{\textcolor{red}{433}}&&\bullet\ar@{-}[ldd]\ar@{-}[rdd]
&&&&\bullet\ar@{-}[ldd]\ar@{-}[rdd]
&&{\textcolor{red}{169}}&&\bullet\ar@{-}[ldd]\ar@{-}[rdd]\\
&&&&{\textcolor{blue}{130}}&&&&&&&&{\textcolor{blue}{894}}&&&&&&&&{\textcolor{blue}{1908}}&&
&&&&&&{\textcolor{blue}{580}}&&&&&\\
&&{\textcolor{red}{\scriptstyle89}}
&&&&{\textcolor{red}{\scriptstyle1325}}
&&&&{\textcolor{red}{\scriptstyle7561}}
&&&&{\textcolor{red}{\scriptstyle2897}}
&&&&{\textcolor{red}{\scriptstyle6466}}
&&&&{\textcolor{red}{\scriptstyle37666}}
&&&&{\textcolor{red}{\scriptstyle14701}}
&&&&{\textcolor{red}{\scriptstyle985}}&&\\
&&{\textcolor{blue}{\scriptstyle420}}
&&&&{\textcolor{blue}{\scriptstyle8503}}
&&&&{\textcolor{blue}{\scriptstyle54450}}
&&&&{\textcolor{blue}{\scriptstyle18222}}
&&&&{\textcolor{blue}{\scriptstyle39366}}
&&&&{\textcolor{blue}{\scriptstyle256050}}
&&&&{\textcolor{blue}{\scriptstyle85610}}
&&&&{\textcolor{blue}{\scriptstyle4077}}
\\
&&&&\ldots&&&&&&&&&&&&\ldots&&&&&&&&&&&&\ldots
}
$$
\textcolor{red}{Markoff's tree}\; /\; \textcolor{blue}{Shadow}

\medskip

Numeric computations confirm that the first hundred terms of the shadow part of this tree are positive numbers.
The ``shadow part'' of the (left) Fibonacci branch branch starts with $0,2,10,38,130,420, 1308, 3970,\ldots$ 
turns out to coincide with Sequence \seqnum{A281199} of the OEIS; see~\cite{OEIS}.
No other branch of this tree has appeared in the OEIS so far.

\subsection{The tree of the vertex $(1:0:2)$}

The second vertex of the quadrilateral corresponds to the following tree.

$$
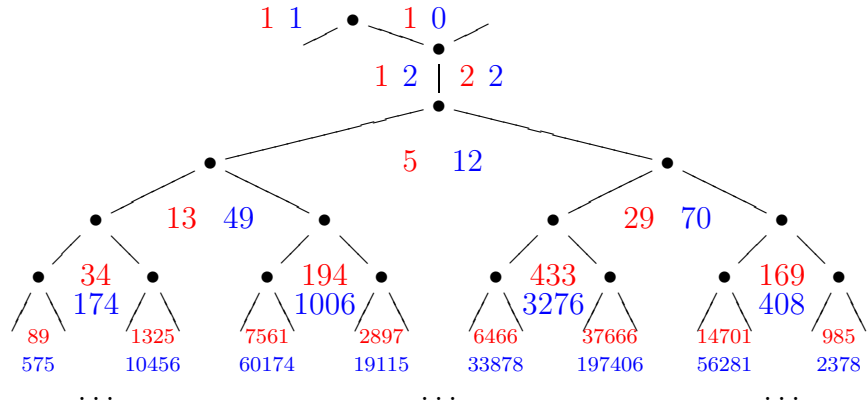
\xymatrix @!0 @R=0.38cm @C=0.38cm
{&&&&&&&&&&&&\\
&&&&&&&&&&{\textcolor{red}{1}}&{\textcolor{blue}{1}}
&&\bullet\ar@{-}[lld]\ar@{-}[lu]&&{\textcolor{red}{1}}&{\textcolor{blue}{0}}&&\\
&&&&&&&&&&&&&&&&\bullet\ar@{-}[dd]\ar@{-}[rru]\ar@{-}[lllu]\\
&&&&&&&&&&&&&&{\textcolor{red}{1}}&{\textcolor{blue}{2}}&&{\textcolor{red}{2}}&{\textcolor{blue}{2}}\\
&&&&&&&&&&&&&&&&\bullet\ar@{-}[lllllllldd]\ar@{-}[rrrrrrrrdd]&&&&&&&&\\
&&&&&&&&&&&&&&&&\\
&&&&&&&&\bullet\ar@{-}[lllldd]\ar@{-}[rrrrdd]
&&&&&&&{\textcolor{red}{5}}&&{\textcolor{blue}{12}}
&&&&&&&\bullet\ar@{-}[lllldd]\ar@{-}[rrrrdd]&&&\\
&&&&&&&&
&&&&&&&&&&&&&&&&\\
&&&&\bullet\ar@{-}[lldd]\ar@{-}[rrdd]
&&&{\textcolor{red}{13}}&&{\textcolor{blue}{49}}&&&\bullet\ar@{-}[lldd]\ar@{-}[rrdd]
&&&&&&&&\bullet\ar@{-}[lldd]\ar@{-}[rrdd]
&&&{\textcolor{red}{29}}&&{\textcolor{blue}{70}}&&&\bullet\ar@{-}[lldd]\ar@{-}[rrdd]\\
&&&&&&&&&&&&&&&&&&&&&&&&&&&&\\
&&\bullet\ar@{-}[ldd]\ar@{-}[rdd]
&&{\textcolor{red}{34}}&&\bullet\ar@{-}[ldd]\ar@{-}[rdd]
&&&&\bullet\ar@{-}[ldd]\ar@{-}[rdd]
&&{\textcolor{red}{194}}&&\bullet\ar@{-}[ldd]\ar@{-}[rdd]
&&&&\bullet\ar@{-}[ldd]\ar@{-}[rdd]
&&{\textcolor{red}{433}}&&\bullet\ar@{-}[ldd]\ar@{-}[rdd]
&&&&\bullet\ar@{-}[ldd]\ar@{-}[rdd]
&&{\textcolor{red}{169}}&&\bullet\ar@{-}[ldd]\ar@{-}[rdd]\\
&&&&{\textcolor{blue}{174}}&&&&&&&&{\textcolor{blue}{1006}}&&&&&&&&{\textcolor{blue}{3276}}&&
&&&&&&{\textcolor{blue}{408}}&&&&&\\
&&{\textcolor{red}{\scriptstyle89}}
&&&&{\textcolor{red}{\scriptstyle1325}}
&&&&{\textcolor{red}{\scriptstyle7561}}
&&&&{\textcolor{red}{\scriptstyle2897}}
&&&&{\textcolor{red}{\scriptstyle6466}}
&&&&{\textcolor{red}{\scriptstyle37666}}
&&&&{\textcolor{red}{\scriptstyle14701}}
&&&&{\textcolor{red}{\scriptstyle985}}&&\\
&&{\textcolor{blue}{\scriptstyle575}}
&&&&{\textcolor{blue}{\scriptstyle10456}}
&&&&{\textcolor{blue}{\scriptstyle60174}}
&&&&{\textcolor{blue}{\scriptstyle19115}}
&&&&{\textcolor{blue}{\scriptstyle33878}}
&&&&{\textcolor{blue}{\scriptstyle197406}}
&&&&{\textcolor{blue}{\scriptstyle56281}}
&&&&{\textcolor{blue}{\scriptstyle2378}}
\\
&&&&\ldots&&&&&&&&&&&&\ldots&&&&&&&&&&&&\ldots
}
$$
Once again, numeric computations confirm the positivity of the shadow part.

Surprisingly, the shadow of the (right)  branch of odd Pell numbers starting with 
$$0, 2, 12, 70, 408, 2378, 13860, 80782,\ldots$$
is nothing else but the sequence of even Pell numbers $P_{2n}$; see Sequence \seqnum{A001542} of~\cite{OEIS}.
No other branch of this tree has been recognized so far.

\subsection{The case of the vertex $(1:1:1)$: double Markoff tree}

It turns out that when the initial conditions for $(\a,\beta,\g)$ are $(1,1,1)$,
as the same as in the classical Markoff tree,
the shadow part doubles the classical one.

\begin{lem}
\label{DoubLem}
The shadow part of the tree with the root
$$
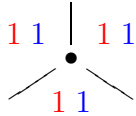
\xymatrix @!0 @R=0.3cm @C=0.3cm
{
&&&\ar@{-}[ddd]\\
\\
&\textcolor{red}{1}\;\textcolor{blue}{1}&&&&\textcolor{red}{1}\;\textcolor{blue}{1}\\
&&&\bullet\ar@{-}[llldd]\ar@{-}[rrrdd]\\
\\
&&&\textcolor{red}{1}\;\textcolor{blue}{1}&&&
}
$$
coincides with the classical Markoff tree.
\end{lem}

\begin{proof}
Given a triple $(A,B,C)=(a+\alpha\epsilon, b+\beta\epsilon, c+\gamma\epsilon)$, then
after the mutation at $A$, the shadow part begins.
$$
\alpha'=\dfrac{2b\beta+2c\gamma-a'\alpha}{a},
$$
see~\cite{Ovs,OB}.
Therefore, when $\a=a,\beta=b,\g=c$, one has
$$
\alpha'=
\dfrac{2b^2+2c^2}{a}-a'=
\dfrac{b^2+c^2}{a},
$$
thanks to~\eqref{muta}.
Hence the lemma.
\end{proof}

Although this vertex does not give an interesting tree, this is the only case for which
positivity is proved.

\subsection{The tree of the vertex $(0:2:1)$}

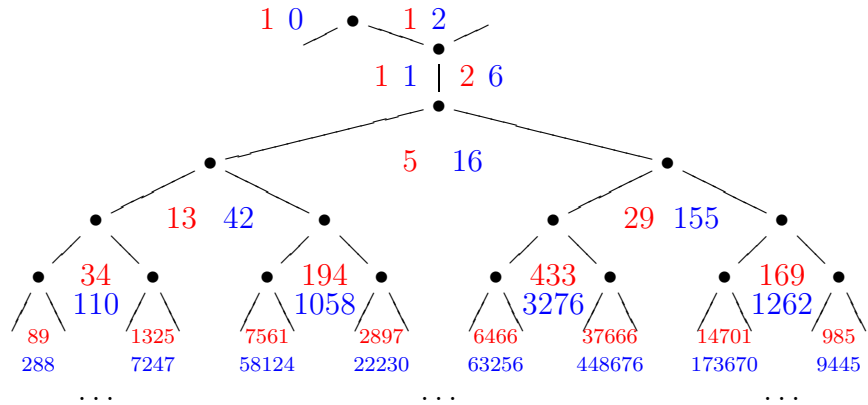
\begin{center}
$$
\xymatrix @!0 @R=0.38cm @C=0.38cm
{&&&&&&&&&&&&\\
&&&&&&&&&&{\textcolor{red}{1}}&{\textcolor{blue}{0}}
&&\bullet\ar@{-}[lld]\ar@{-}[lu]&&{\textcolor{red}{1}}&{\textcolor{blue}{2}}&&\\
&&&&&&&&&&&&&&&&\bullet\ar@{-}[dd]\ar@{-}[rru]\ar@{-}[lllu]\\
&&&&&&&&&&&&&&{\textcolor{red}{1}}&{\textcolor{blue}{1}}&&{\textcolor{red}{2}}&{\textcolor{blue}{6}}\\
&&&&&&&&&&&&&&&&\bullet\ar@{-}[lllllllldd]\ar@{-}[rrrrrrrrdd]&&&&&&&&\\
&&&&&&&&&&&&&&&&\\
&&&&&&&&\bullet\ar@{-}[lllldd]\ar@{-}[rrrrdd]
&&&&&&&{\textcolor{red}{5}}&&{\textcolor{blue}{16}}
&&&&&&&\bullet\ar@{-}[lllldd]\ar@{-}[rrrrdd]&&&\\
&&&&&&&&
&&&&&&&&&&&&&&&&\\
&&&&\bullet\ar@{-}[lldd]\ar@{-}[rrdd]
&&&{\textcolor{red}{13}}&&{\textcolor{blue}{42}}&&&\bullet\ar@{-}[lldd]\ar@{-}[rrdd]
&&&&&&&&\bullet\ar@{-}[lldd]\ar@{-}[rrdd]
&&&{\textcolor{red}{29}}&&{\textcolor{blue}{155}}&&&\bullet\ar@{-}[lldd]\ar@{-}[rrdd]\\
&&&&&&&&&&&&&&&&&&&&&&&&&&&&\\
&&\bullet\ar@{-}[ldd]\ar@{-}[rdd]
&&{\textcolor{red}{34}}&&\bullet\ar@{-}[ldd]\ar@{-}[rdd]
&&&&\bullet\ar@{-}[ldd]\ar@{-}[rdd]
&&{\textcolor{red}{194}}&&\bullet\ar@{-}[ldd]\ar@{-}[rdd]
&&&&\bullet\ar@{-}[ldd]\ar@{-}[rdd]
&&{\textcolor{red}{433}}&&\bullet\ar@{-}[ldd]\ar@{-}[rdd]
&&&&\bullet\ar@{-}[ldd]\ar@{-}[rdd]
&&{\textcolor{red}{169}}&&\bullet\ar@{-}[ldd]\ar@{-}[rdd]\\
&&&&{\textcolor{blue}{110}}&&&&&&&&{\textcolor{blue}{1058}}&&&&&&&&{\textcolor{blue}{3276}}&&
&&&&&&{\textcolor{blue}{1262}}&&&&&\\
&&{\textcolor{red}{\scriptstyle89}}
&&&&{\textcolor{red}{\scriptstyle1325}}
&&&&{\textcolor{red}{\scriptstyle7561}}
&&&&{\textcolor{red}{\scriptstyle2897}}
&&&&{\textcolor{red}{\scriptstyle6466}}
&&&&{\textcolor{red}{\scriptstyle37666}}
&&&&{\textcolor{red}{\scriptstyle14701}}
&&&&{\textcolor{red}{\scriptstyle985}}&&\\
&&{\textcolor{blue}{\scriptstyle288}}
&&&&{\textcolor{blue}{\scriptstyle7247}}
&&&&{\textcolor{blue}{\scriptstyle58124}}
&&&&{\textcolor{blue}{\scriptstyle22230}}
&&&&{\textcolor{blue}{\scriptstyle63256}}
&&&&{\textcolor{blue}{\scriptstyle448676}}
&&&&{\textcolor{blue}{\scriptstyle173670}}
&&&&{\textcolor{blue}{\scriptstyle9445}}
\\
&&&&\ldots&&&&&&&&&&&&\ldots&&&&&&&&&&&&\ldots
}
$$
\end{center}
The positivity of the shadow part is also confirmed numerically.

This time, one can see the left branch 
starting with $2,6,16,42,110,288,754, 1974,\ldots$ corresponds to Sequence \seqnum{A025169}
consisting of the numbers $2F_{2n+2}$.
No other branch of this tree is recognized so far.

\section{Code and numeric evidence for the conjecture}\label{Python}

In this section, I will explain how the domain of possible solutions has been reduced and demonstrate that deviating from the point $(1,1,1)$ in ``wrong'' directions
leads to negative numbers in the shadow Markoff tree.
Similar computations are also working for the other vertices of the quadrilateral of Conjecture~\ref{MainConj}.

I can restrict the possibilities for the position of the convex set because of
an observation.
The condition $\alpha\ge0$ can be written $\dfrac{-\a a'+2b\beta+2c\gamma}{a}\ge0$.
Following the left branch of the Markoff tree where one element of the triple
does not move (let us say $c+\e \gamma= 1 +\e$), one can transform the condition to be
$$
\dfrac{2b\beta +2}{a' a}\le \alpha.
$$
This means $(\alpha,\beta)$ has to stay left of a certain linear function, cutting the plan
in half. Tracing those functions while browsing the left branch of the tree gives lots
of conditions, and the following graph.

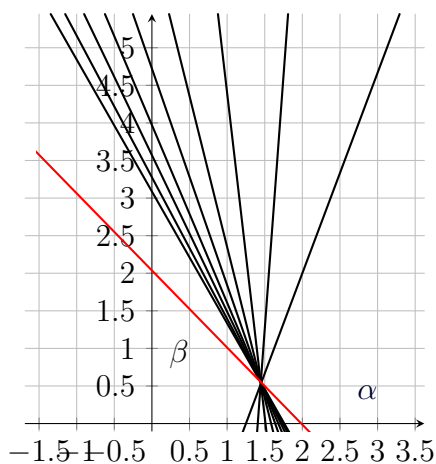
\begin{center}
\begin{tikzpicture}[line cap=round,line join=round,>=triangle 45,x=1cm,y=1cm]
\begin{axis}[
x=1cm,y=1cm,
axis lines=middle,
ymajorgrids=true,
xmajorgrids=true,
ymin=-0.1,
ymax=5.44822254231216,
xmin=-1.6844120637227598,
xmax=3.621939686427615,
ytick={-1.5,-1,...,5},
xtick={-1.5,-1,...,3.5},
]
\clip(-1.5370441957610437,-1.6844120637227598) rectangle (5.44822254231216,5.621939686427615);
\draw[line width=0.8pt,smooth,samples=100,domain=-1.5370441957610437:5.44822254231216] plot({0.375*(\x)+1.25},\x);
\draw[line width=0.8pt,smooth,samples=100,domain=-1.5370441957610437:5.44822254231216] plot({0.074*(\x)+1.407},\x);
\draw[line width=0.8pt,smooth,samples=100,domain=-1.5370441957610437:5.44822254231216] plot({0-0.116*(\x)+1.51},\x);
\draw[line width=0.8pt,smooth,samples=100,domain=-1.5370441957610437:5.44822254231216] plot({0-0.249*(\x)+1.5849},\x);
\draw[line width=0.8pt,smooth,samples=100,domain=-1.5370441957610437:5.44822254231216] plot({0-0.347*(\x)+1.639},\x);
\draw[line width=0.8pt,smooth,samples=100,domain=-1.5370441957610437:5.44822254231216] plot({0-0.423*(\x)+1.68},\x);
\draw[line width=0.8pt,smooth,samples=100,domain=-1.5370441957610437:5.44822254231216] plot({0-0.48*(\x)+1.714},\x);
\draw[line width=0.8pt,smooth,samples=100,domain=-1.5370441957610437:5.44822254231216] plot({0-0.531*(\x)+1.74},\x);
\draw[line width=0.8pt,smooth,samples=100,domain=-1.5370441957610437:5.44822254231216] plot({0-0.57*(\x)+1.76},\x);
\draw[line width=0.8pt,color=ffqqqq,smooth,samples=100,domain=-1.5370441957610437:5.44822254231216] plot({0-0.975*(\x)+1.986},\x);
\draw [color=qqqqtt](2.587933756476852,0.6538160252438828) node[anchor=north west] {$\alpha$};
\draw (0.087933756476852,1.2538160252438828) node[anchor=north west] {$\beta$};
\end{axis}
\end{tikzpicture}
\end{center}
In black, you can notice the first restricting lines. 
The red line is the limit that appears when drawing hundreds of those.
The solutions had to be left of all these lines, it was then necessary to refine the domain.

Let us show with computations that deviating from the point $(1,1,1)$ in ``wrong'' directions leads to negative numbers in the shadow Markoff tree.\\

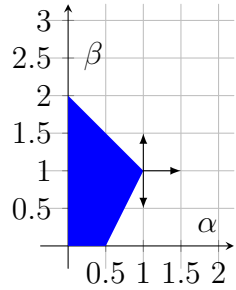
\begin{center}
\definecolor{zzttqq}{rgb}{0,0,1}
\begin{tikzpicture}[line cap=round,line join=round,>=triangle 45,x=1cm,y=1cm]
\begin{axis}[
x=1cm,y=1cm,
axis lines=middle,
ymajorgrids=true,
xmajorgrids=true,
ymin=-0.29868247168517575,
ymax=3.2156316474709685,
xmin=-0.36602351946208256,
xmax=2.204001897585711,
ytick={0,0.5,...,3},
xtick={0,0.5,...,2},]
\clip(-0.36602351946208256,-0.29868247168517575) rectangle (4.204001897585711,3.2156316474709685);
\fill[line width=2pt,color=zzttqq,fill=zzttqq,fill opacity=0.50000000149011612] (0,0) -- (0.5,0) -- (1,1) -- (0,2) -- cycle;
\draw [->,>=latex] (1.,1.) -- (1.5,1.);
\draw [->,>=latex] (1.,1.) -- (1.,1.5);
\draw [->,>=latex] (1.,1.) -- (1.,0.5);
\draw (0.06610799040603345,2.8605427849359513) node[anchor=north west] {$\beta$};
\draw (1.5738817601573287,0.50538827311252061) node[anchor=north west] {$\alpha$};
\end{axis}
\end{tikzpicture}
\end{center}

I present a function path(alpha,beta,gamma,C), where (alpha,beta,gamma) is the shadow part of the root of the tree, and C is a list composed of l (for left) and r (for right). It gives the list of triple encountered in the tree when following this path after the 2 mandatory beginning mutations (to the right, then to the left).\\
The language used is sage.

\begin{python}

   def path(alpha,beta,gamma,C):#C is the path to take, in a list.
    (a,b,c)=(1,1,1)
    C='r'+'l'+C		#l for left, r for right.
    L=[[a,alpha,b,beta,c,gamma]]#list of dual numbers in the tree
    for i in C:
        if i=='l':#mutation on the second element of the triple:beta
            (a,b,c,alpha,beta,gamma)=(a,c, (a^2+c^2)/b , alpha, gamma, (-((a^2+c^2)/b)*beta +2*a*alpha+2*c*gamma)/b )
            L.append([a,alpha,b,beta,c,gamma])
        if i=='r':#mutation on the first element of the triple: alpha
            (a,b,c,alpha,beta,gamma)=(c,b, (b^2+c^2)/a , gamma, beta, (-((b^2+c^2)/a)*alpha +2*b*beta+2*c*gamma)/a )
            L.append([a,alpha,b,beta,c,gamma])
    return(L)					
					\end{python}

The answer is in the form of a list, each element of the list is presented as 

\begin{center}
[a,alpha,b,beta,c,gamma].
\end{center}

\subsection{Top side}

$path(1,0.9,1,'rlrlrlr')$ gives a negative number in the last triple\\
 $\gamma\approx -6.98e33$.
It shows that it is not possible to go to the top of the (1,1,1) point.

One can access the last triple directly using 
\begin{python}
path(1,0.9,1,'rlrlrlr')[9]
\end{python}
or only the value of the negative $\gamma$ with
\begin{python}
path(1,0.9,1,'rlrlrlr')[9][5]
\end{python}

\subsection{Bottom side}

$path(1,1.1,1,'llllrrrrrrrrr')$ (4 times ``l" and 9 times ``r") gives a negative number in the last triple $\gamma\approx -1.38e23$.\\
It shows that it is not possible to go to under the (1,1,1) point.

It is possible to access the last triple directly using 
\begin{python}
path(1,1.1,1,'llllrrrrrrrrr')[15]
\end{python}
or only the value of the negative $\gamma$ with
\begin{python}
path(1,1.1,1,'llllrrrrrrrrr')[15][5]
\end{python}

\subsection{Up side}

$path(1.1,1,1,'llllrr')$ gives a negative number in the last triple $\gamma\approx -77761.8$.\\
It shows that it is not possible to go to the right of the (1,1,1) point.

You can access the last triple directly using 
\begin{python}
path(1.1,1,1,'llllrr')[8]
\end{python}
or only the value of the negative $\gamma$ with
\begin{python}
path(1.1,1,1,'llllrr')[8][5]
\end{python}

I also wrote a function that constructs a full tree of height $n$.
The function is 
\begin{center}
shadow(alpha,beta,gamma,n), 
\end{center}
where (alpha,beta,gamma) is the root of the tree, and n is the height of the shadow tree to build. It displays the shadow tree in the form of a list, beginning after the 2 mandatory mutations.
The tree is defined as a list of 3 elements: [[a,alpha,b,beta,c,gamma],left son, right son]. The left and right sons are also trees of the same form, or are empty lists: [ ]\\
\begin{python}
[[1, 0, 1, 0, 1, 1],
 [1, 1, 1, 0, 2, 2],
 [1, 1, 2, 2, 5, 10],
 [1, 1, 5, 10, 13, 38],
 [1, 1, 13, 38, 34, 130],
 [1, 1, 34, 130, 89, 420]]	
					\end{python}

For example shadow(0,0,1,3) gives

\begin{python} 
[[1, 1, 2, 2, 5, 10],
 [[1, 1, 5, 10, 13, 38],
  [[1, 1, 13, 38, 34, 130], [], []],
  [[13, 38, 5, 10, 194, 894], [], []]],
 [[5, 10, 2, 2, 29, 79],
  [[5, 10, 29, 79, 433, 1908], [], []],
  [[29, 79, 2, 2, 169, 580], [], []]]]
					\end{python}

I use multiple steps to generate the tree. First I generate a binary tree of height $n$, then I fill it using the relations for $(a',\alpha')$, beginning after the 2 mandatory mutations.

\begin{python} 
def binary(alpha,beta,gamma,n):#gives a binary tree of hight n
    if n==0:
        return([])
    if n>0:
        return([[1,alpha,2,beta,5,gamma],binary(alpha,beta,gamma,n-1),binary(alpha,beta,gamma,n-1)])

def genere(arbre):#Put the shadow Markoff numbers in the tree
    if not(arbre[1]==[]):#order: a, alpha, b, beta, c, gamma
        arbre[1][0][0]=arbre[0][0]
        arbre[1][0][1]=arbre[0][1]
        arbre[1][0][2]=arbre[0][4]
        arbre[1][0][3]=arbre[0][5]      
        arbre[1][0][4]=(arbre[0][0]^2+arbre[0][4]^2)/(arbre[0][2])
        arbre[1][0][5]=(-arbre[1][0][4]*arbre[0][3]+2*arbre[0][0]*arbre[0][1]+2*arbre[0][4]*arbre[0][5])/(arbre[0][2])
        genere(arbre[1])
    if not(arbre[2]==[]):        
        arbre[2][0][0]=arbre[0][4]
        arbre[2][0][1]=arbre[0][5]
        arbre[2][0][2]=arbre[0][2]
        arbre[2][0][3]=arbre[0][3]      
        arbre[2][0][4]=(arbre[0][2]^2+arbre[0][4]^2)/(arbre[0][0])
        arbre[2][0][5]=(-arbre[2][0][4]*arbre[0][1]+2*arbre[0][2]*arbre[0][3]+2*arbre[0][4]*arbre[0][5])/(arbre[0][0])
        genere(arbre[2])

def shadow(alpha,beta,gamma,n):
    (a,b,c)=(1,1,1)
    (a,b,c,alpha,beta,gamma)=(c,b, (b^2+c^2)/a , gamma, beta, (-((b^2+c^2)/a)*alpha +2*b*beta+2*c*gamma)/a )
    (a,b,c,alpha,beta,gamma)=(a,c, (a^2+c^2)/b , alpha, gamma, (-((a^2+c^2)/b)*beta +2*a*alpha+2*c*gamma)/b )
    A=binary(alpha,beta,gamma,n)
    genere(A)
    return(A)

					\end{python}

\medskip
 \noindent
\section{Acknowledgments}
I am grateful to Valentin Ovsienko for constant help and many enlightening discussions.

2020 {\it Mathematics Subject Classification}:
Primary 11D25; Secondary 15A75, 13F60.

\noindent \emph{Keywords: } Markoff tree, dual number, positivity

\bigskip
\hrule
\bigskip

\noindent (Concerned with sequence
\seqnum{A000045}, \seqnum{A000129}, \seqnum{A001542}, \seqnum{A025169}, \seqnum{A002559}, \seqnum{A238846}, and \seqnum{A281199}.)

\bigskip
\hrule
\bigskip

\vspace*{+.1in}
\noindent
Received 
revised version received  
Published in {\it Journal of Integer Sequences},
\bigskip
\hrule
\bigskip

\noindent
Return to 
\href{https://cs.uwaterloo.ca/journals/JIS/}{Journal of Integer Sequences home page} 
\vskip .1in

\end{document}